# Multiscale seamless-domain method for nonperiodic fields: linear heat conduction analysis

Yoshiro Suzuki

Tokyo Institute of Technology, Department of Mechanical Engineering,
2-12-1 Ookayama, Meguro-ku, Tokyo 152-8552, Japan
e-mail address: ysuzuki@ginza.mes.titech.ac.jp

**Abstract**

A multiscale numerical solver called the seamless-domain method (SDM) is used in linear heat conduction analysis of nonperiodic simulated fields. The practical feasibility of the SDM has been verified for use with periodic fields but has not previously been verified for use with nonperiodic fields. In this paper, we illustrate the mathematical framework of the SDM and the associated error factors in detail. We then analyze a homogeneous temperature field using the SDM, the standard finite difference method, and the conventional domain decomposition method (DDM) to compare the convergence properties of these methods. In addition, to compare their computational accuracies and time requirements, we also simulated a nonperiodic temperature field with a nonuniform thermal conductivity distribution using the three methods. The accuracy of the SDM is very high and is approximately equivalent to that of the DDM. The mean temperature error is less than 0.02% of the maximum temperature in the simulated field. The total central processing unit (CPU) times required for the analyses using the most efficient SDM model and the DDM model represent 13% and 17% of that of the finite difference model, respectively.

## 1. Theory of the SDM

The seamless-domain method (SDM) is a multiscale numerical technique composed of macroscopic global analysis and microscopic local analysis methods. If we add a mesoscopic analysis between the global and local simulations, the SDM then becomes a three-scale analysis method. In this study, we perform two-scale SDM analyses of linear stationary temperature fields.

The total simulated field is called the global domain, and is denoted by $G$. $G$ has no meshes, elements, grids, or cells. The SDM is therefore categorized as a mesh-free method.

As shown in Fig. 1(a), $G$ is represented using coarse-grained points (CPs). These CPs are used to spatially discretize $G$. The distribution of the material properties in $G$ is arbitrary. Figure 1 shows an example of $G$ with a nonuniform and nonperiodic material property profile. Even if $G$ is heterogeneous and has microscopic constituents, we do not need to model them separately. Therefore, the number of microscopic inclusions that can be located between adjacent CPs is arbitrary.

As shown in Fig. 1, $G$ is divided into multiple small domains that are called local domains ($L$). Each local domain has a CP at its center and other neighboring CPs. Local analysis of $L$ provides the relationships among the CPs (i.e., it provides the SDM's discretized equation, which is described in Subsection 2.1). The relationship estimates the dependent-variable value(s) at the center CP of each domain with reference to the variables at the neighboring CPs. By formulating and solving these relational expressions for all CPs in $G$, we can obtain all the variables for the domain.

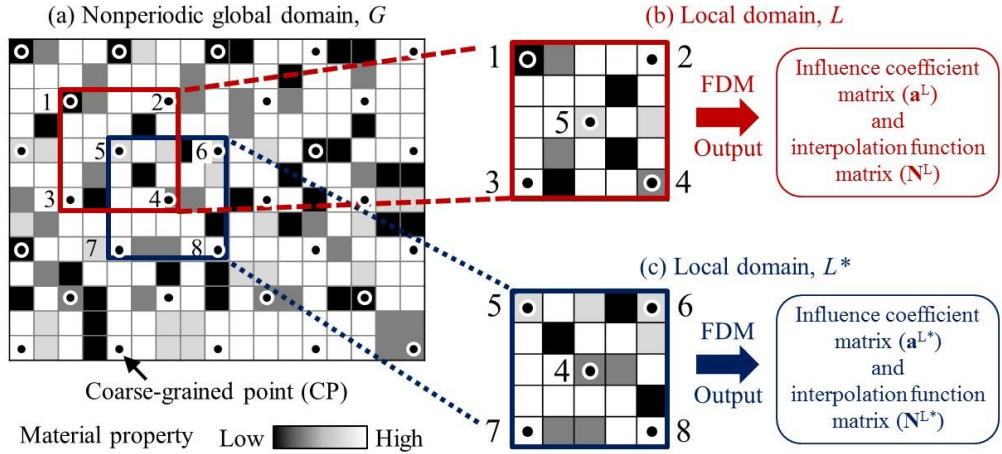

Fig. 1. (a) Nonperiodic and heterogeneous global domain; (b) the SDM's local domain $L$; and (c) a local domain adjacent to $L$.

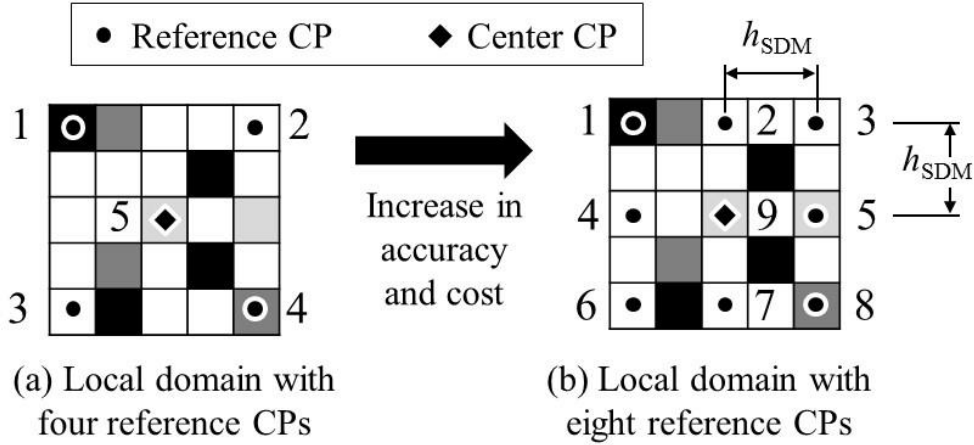

Fig. 2. (a) Local domain with four reference CPs, and (b) local domain with eight reference CPs.

The local analysis can be performed using conventional numerical schemes such as the finite difference method (FDM), the finite element method, and the finite volume method. The FDM is used in the work in this manuscript.

The SDM has previously only been applied to periodic fields that consist of unit cells with identical structures. It is demonstrated here that the SDM can quickly reproduce solutions that are as accurate as those obtained using fine-grained FEM and FDM models at low cost for both linear problems involving periodic fields (including stationary heat conduction [1–3], nonstationary heat conduction [4], elasticity [2,5]) and nonlinear problems (e.g., stationary heat conduction [6]).

In contrast, this study applies the SDM to a nonperiodic field (as shown in Fig. 1(a)) for the first time. We demonstrate how to analyze a field using the SDM in cases where the distribution of the material properties in the field is both nonuniform and nonperiodic.

In Section 2, we describe the theory of the SDM in detail and demonstrate discretization of the governing equation(s). Additionally, we explain the errors that are generated by the discretization process.

Section 3 provides the computational implementation when a standard FDM solver is used to perform the local analysis for the SDM.

Section 4 describes a convergence study for a linear stationary temperature problem. The target field is a square with an isotropic and uniform thermal conductivity distribution. We test the SDM models along with FDM models and domain decomposition method (DDM) models. Like the SDM, the DDM divides the entire field into multiple local fields to reduce the total computational time. This problem can be solved completely and we can thus obtain the true temperature. We use temperature differences when compared with the true solution as an accuracy index to compare the convergence properties of the three methods.

Section 5 gives an example of a nonperiodic field problem. We use the three methods (i.e., SDM, FDM, and DDM) to analyze a field with an isotropic, nonuniform, and nonperiodic thermal conductivity distribution. We cannot obtain the true solution to this problem. Therefore, by comparison with the FDM solution, we investigate the calculation accuracies of both the SDM and the DDM.

## 2 Theory of the SDM
## 2.1 Outline of the SDM

Section 2 illustrates the theory of the SDM using the example of the following elliptic partial differential equation. The theory is generally applicable to any other linear problem.

$$\frac{\partial}{\partial x}\left(k\frac{\partial u}{\partial x}\right) + \frac{\partial}{\partial y}\left(k\frac{\partial u}{\partial y}\right) = -f \text{ in } G, \quad (1)$$

where $u = u(\mathbf{x})$ is the dependent-variable value (i.e., temperature in a heat conduction problem) at point $\mathbf{x} = (x,y)^T$, $k = k(\mathbf{x})$ is the conductivity (specifically, the thermal conductivity here), and $f = f(\mathbf{x})$ is the source value (i.e., the calorific value). Note that $\#^T$ represents the transpose of $\#$.

As shown in Fig. 1(a), the complete simulated field is called the global domain and is denoted by $G$. In the SDM, we arrange the CPs in $G$ and thus express $G$ using only the information on the CPs, i.e., $G$ is spatially discretized using the CPs.

Subsequently, we divide $G$ into multiple local domains (denoted by $L$). Figure 2(a) shows an example of $L$. Each local domain has a CP at its center that is surrounded by other CPs. In general, the number of CPs contained in $G$ is equal to the number of local domains in $G$.

The $L$ that is depicted in Fig. 2(a) has CP 5 (rhombus) at its center and four reference CPs (CPs 1–4, circles) arranged around CP 5. The dependent-variable value of CP $i$ is denoted by $u_{\text{CP}i}$:

$$u_{\text{CP}i} = u(\mathbf{x}_{\text{CP}i}). \quad (2)$$

When the governing equation (Eq. (1)) is linear, the relationship between CP 5 and CPs 1–4 (i.e., the discretized equation for $u$) can be expressed in the following form.

$$u_{CP5} = a_1^L u_{CP1} + a_2^L u_{CP2} + a_3^L u_{CP3} + a_4^L u_{CP4} + F^L(\mathbf{x}_{CP5})$$
$$= \mathbf{a}^L \mathbf{u}^L + F^L(\mathbf{x}_{CP5}), \quad (3)$$

where
$$\mathbf{u}^L = (u_{CP1} \cdots u_{CP4})^T. \quad (4)$$

$F^L(\mathbf{x}_{CP5})$ is a constant term and has different values that depend on the source distribution, $f(\mathbf{x})$, in $L$. $a_i^L$ represents the weighting factor for CP $i$ (and is referred to here as the influence coefficient). In the SDM, the one-by-four matrix shown below, which includes influence coefficients for CPs 1–4,
$$\mathbf{a}^L = (a_1^L \cdots a_4^L), \quad (5)$$
is called the influence coefficient matrix for $L$.

$\mathbf{a}^L$ is constructed from the results of the local analysis of $L$. A detailed illustration of the derivation of $\mathbf{a}^L$ is provided in Subsections 2.2 and 3.1.

We formulate an equation in the form of Eq. (3) for each local domain and then construct simultaneous equations. The unknown variable values in all CPs can subsequently be determined by solving these simultaneous equations.

The most notable characteristic of the SDM is that all adjacent local domains are partially superposed with each other and thus share some of their CPs. For example, Fig. 1(b) and (c) show that $L$ and $L^*$ (which denotes the local domain next to $L$) overlap each other's quarter regions and share CPs 4 and 5. The SDM guarantees that the variables for the shared CPs ($u_{CP4}, u_{CP5}$) contained in $L$ are exactly equal to those contained in $L^*$. The shared CPs therefore attempt to equalize the variable distributions of the overlapped region in $L$ and the corresponding region in $L^*$.

For example, if we increase the number of reference CPs from four to eight, as shown in Fig. 2(b), the number of shared CPs then increases to four. This increase in the number of shared CPs would also enhance the variable continuity between $L$ and $L^*$.

Because the number of reference CPs is arbitrary, we can increase the numbers of reference CPs according to the computational precision required. However, any increase in the number of reference CPs also leads to an increase in the computation time.

## 2.2 Local analysis and calculation errors of the SDM

### 2.2.1 Influence coefficient matrix $\mathbf{a}^L$ and interpolation function matrix $\mathbf{N}^L(\mathbf{x})$

We use $L$ as an example in this subsection. $L$ has four CPs (CPs 1–4) on its boundary, $\partial L$.

First, we give the interpolation functions for $L$ in the form $\mathbf{N}^L(\mathbf{x})$. $\mathbf{N}^L(\mathbf{x})$ is a function matrix that is used to estimate the variables of an arbitrary point $\mathbf{x}$ in $L$ (i.e., $u(\mathbf{x})$) with reference to the variables of CPs 1–4, where $\mathbf{u}^L = (u_{CP1} \cdots u_{CP4})^T$.
$$u(\mathbf{x}) = \mathbf{N}^L(\mathbf{x})\mathbf{u}^L + F^L(\mathbf{x}). \quad (6)$$

By substituting $\mathbf{x} = \mathbf{x}_{CP5}$ (i.e., the position of CP 5) into the equation above, it follows that
$$u_{CP5} = u(\mathbf{x}_{CP5}) = \mathbf{N}^L(\mathbf{x}_{CP5})\mathbf{u}^L + F^L(\mathbf{x}_{CP5}). \quad (7)$$

Comparison of Eq. (3) and Eq. (7) allows $\mathbf{a}^L$ to be expressed in the form below:
$$\mathbf{a}^L = \mathbf{N}^L(\mathbf{x}_{CP5}). \quad (8)$$

Therefore, we can easily compute $\mathbf{a}^L$ from $\mathbf{N}^L$.

The derivation of $\mathbf{N}^L$ is described below. Because $L$ is a linear field, the function of $u(\mathbf{x})$ in $L$ is uniquely determined when the following two items are provided:
· the boundary condition of $L$, e.g., $u$ on $\partial L$;
· the source distribution, i.e., $f(\mathbf{x})$ in $L$.
$$u(\mathbf{x}) = \int_{\partial L} b(\mathbf{x}, s) u(s) ds + \int_L c(\mathbf{x}, \mathbf{x}') f(\mathbf{x}') d\mathbf{x}', \quad (9)$$
where $b$ and $c$ are unknown functions. Because $L$ is a two-dimensional field in this case, the first term on the right side is a line integral and the second term is a surface integral. If $L$ is three-dimensional, then the first and second terms are a surface integral and a volume integral, respectively.

The first term on the right side is the line integral of the entire closed curve $\partial L$ (i.e., the boundary of $L$). $u(s)$ is the variable value at point $s$ on $\partial L$, where $s$ is the distance on $\partial L$ from a specific point, and $ds$ is the line element.

Under certain restrictive conditions, exact solutions for $b$ and $c$ can be obtained. For example, if the conductivity is constant and is independent of the position $\mathbf{x}$ (where $k(\mathbf{x}) = const.$), the governing equation (Eq. (1)) then becomes Poisson's equation. In this case, $c$ is expressed using the Green's function, $g$, shown below.
$$c(\mathbf{x}, \mathbf{x}') = g(\mathbf{x}, \mathbf{x}'). \quad (10)$$

However, because we cannot obtain solutions for $b$ and $c$ in many general cases, we compute Eq. (9) approximately using a numerical simulation solver (e.g., the FDM in this paper).

The necessary conditions for exact computation of $u(\mathbf{x})$ are given as follows.
(1) The two integrals on the right side of Eq. (9) can be calculated exactly.
(2) The source distribution in the local domain, i.e., $f(\mathbf{x})$ in $L$, can be obtained.
(3) The boundary condition of $L$, i.e., $u(s)$ on $\partial L$, can be obtained.
(4) Exact solutions for $b$ and $c$ on the right side of Eq. (9) can be obtained.

Even if condition (1) is not satisfied, the integration errors can be suppressed to be sufficiently small as long as accurate numerical integration can be performed. Additionally, $f(\mathbf{x})$ is given (i.e., it is known) in general. Consequently, the main error factors of the SDM are related to (3) and (4). Subsections 2.2.2.1 and 2.2.2.2 describe (3) and (4), respectively.

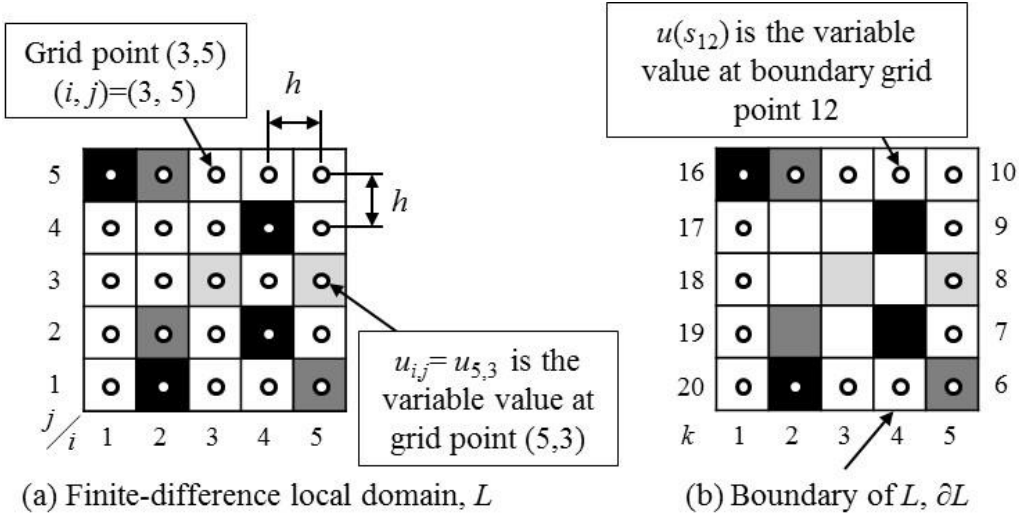

Fig. 3. (a) All grid points in a finite-difference local domain, and (b) boundary grid points in the local domain.

### 2.2.2 Error factors of the SDM
#### 2.2.2.1 Error factor 1: Local boundary condition $u(s)$ on $\partial L$

We consider the case where conditions (2) and (4) are satisfied in this subsection. Therefore, we can obtain exact solutions for $b$, $c$, and $f$ on the right side of Eq. (9).

We define the second term on the right side of Eq. (9) as $F^L(\mathbf{x})$. We can then compute $F^L(\mathbf{x})$ either analytically or numerically.
$$F^L(\mathbf{x}) = \int_L c(\mathbf{x}, \mathbf{x}')f(\mathbf{x}')d\mathbf{x}'. \quad (11)$$

The boundary condition of $L$, denoted by $u(s)$ on $\partial L$, is generally unknown. $\partial L$ is the boundary of $L$. In addition, $L$ is part of the global domain, $G$, and is thus partially superposed on other local domains. $u(s)$ on $\partial L$ is thus determined with a dependence upon the variable profiles in the neighboring local domains. Therefore, we cannot obtain $u(s)$ on $\partial L$ at the local analysis stage.

Let us express $u(s)$ using a mathematical expression by interpolating the variables at CPs 1–4, in the form $\mathbf{u}^L = (u_{\text{CP1}} \cdots u_{\text{CP4}})^T$. $N_i^{\partial L}$ denotes the interpolation function for $u_{\text{CP}i}$. We then obtain
$$u(s) = N_1^{\partial L}(s)u_{\text{CP1}} + N_2^{\partial L}(s)u_{\text{CP2}} + \cdots + N_4^{\partial L}(s)u_{\text{CP4}}$$
$$= \mathbf{N}^{\partial L}(s)\mathbf{u}^L, \quad (12)$$
where
$$\mathbf{N}^{\partial L}(s) = \left(N_1^{\partial L}(s) \cdots N_4^{\partial L}(s)\right) \quad (13)$$
is the interpolation function matrix for $\partial L$. Equation (12) expresses the unknown boundary condition, $u(s)$, using the unknown vector, $\mathbf{u}^L = (u_{\text{CP1}} \cdots u_{\text{CP4}})^T$.

By substituting Eq. (12) into the first term on the right side of Eq. (9), it follows that
$$\int_{\partial L} b(\mathbf{x}, s)u(s)ds = \int_{\partial L} b(\mathbf{x}, s)\mathbf{N}^{\partial L}(s)\mathbf{u}^L ds \quad (14)$$
$$= \int_{\partial L} b(\mathbf{x}, s)\mathbf{N}^{\partial L}(s)ds\, \mathbf{u}^L = \mathbf{N}^L(\mathbf{x})\mathbf{u}^L$$
where
$$\mathbf{N}^L(\mathbf{x}) = \int_{\partial L} b(\mathbf{x}, s)\mathbf{N}^{\partial L}(s)ds \quad (15)$$
is the interpolation function matrix for $L$. Note that $\mathbf{u}^L$ is independent of $s$ and can be separated from the integrand. By substituting Eqs. (11) and (14) into Eq. (9), we obtain
$$u(\mathbf{x}) = \mathbf{N}^L(\mathbf{x})\mathbf{u}^L + F^L(\mathbf{x}). \quad (16)$$
In addition, from Eqs. (8) and (15), we derive the following influence coefficient matrix for CP 5:
$$\mathbf{a}^L = \mathbf{N}^L(\mathbf{x}_{\text{CP5}}) = \int_{\partial L} b(\mathbf{x}_{\text{CP5}}, s)\mathbf{N}^{\partial L}(s)ds. \quad (17)$$
Then,
$$u_{\text{CP5}} = u(\mathbf{x}_{\text{CP5}})$$
$$= \int_{\partial L} b(\mathbf{x}_{\text{CP5}}, s)u(s)ds + \int_L c(\mathbf{x}_{\text{CP5}}, \mathbf{x}')f(\mathbf{x}')d\mathbf{x}' \quad (18)$$
$$= \mathbf{a}^L \mathbf{u}^L + F^L(\mathbf{x}_{\text{CP5}}).$$

It should be noted here that $\mathbf{N}^{\partial L}$ in Eq. (12) is not determined automatically. Additionally, there are many possible ways to construct $\mathbf{N}^{\partial L}$. $\mathbf{N}^{\partial L}$ can be changed to suit the way in which the local analysis is conducted. Therefore, the user of the SDM determines how $\mathbf{N}^{\partial L}$ is constructed, i.e., the user chooses how to conduct the local analysis. Subsections 3.1.1 and 3.1.2 provide detailed descriptions of how to construct $\mathbf{N}^{\partial L}$.

In general, we cannot generate $\mathbf{N}^{\partial L}$ exactly. As shown in Eq. (13), $\mathbf{N}^{\partial L}$ in this subsection is a one-by-four matrix and has four dependent-variable profile modes. If the actual profile on $\partial L$ has more than four modes or if it includes a high-frequency mode that cannot be expressed using $\mathbf{N}^{\partial L}$, $\mathbf{N}^{\partial L}$ then causes an error. When the interpolation functions for $\partial L$ ($\mathbf{N}^{\partial L}$) have an error, the corresponding functions for $L$ ($\mathbf{N}^L$) thus also include an error (see Eq. (15)). In this case, $u_{\text{CP5}}$, which was calculated using Eq. (18), generates an error even if both $b(\mathbf{x}_{\text{CP5}}, s)$ and $\mathbf{u}^L$ are entirely accurate.

$\mathbf{N}^{\partial L}_{\mathrm{TRUE}}$ and $\mathbf{N}^{\partial L}_{\mathrm{ERROR}}$ denote the accurate interpolation matrix for $\partial L$ and its error matrix, respectively. The $\mathbf{N}^{\partial L}$ that is generated in the local analysis can be expressed in the form
$$\mathbf{N}^{\partial L} = \mathbf{N}^{\partial L}_{\mathrm{TRUE}} + \mathbf{N}^{\partial L}_{\mathrm{ERROR}}. \quad (19)$$
The effect of $\mathbf{N}^{\partial L}_{\mathrm{ERROR}}$ on the accuracy of $u_{\mathrm{CP5}}$ is described in Subsection 2.2.2.3. Basically, if the number of entries in $\mathbf{N}^{\partial L}$ (i.e., the number of CPs in $\partial L$) increases, this means that we can reduce $\mathbf{N}^{\partial L}_{\mathrm{ERROR}}$.

### 2.2.2.2 Error factor 2: Precision of numerical solver used in the local analysis

As stated in Subsection 2.2.1, we cannot derive the functions $b$ and $c$ on the right side of Eq. (9) exactly. We therefore need to prepare approximation functions for $b$ and $c$ using a simulation solver to compute Eq. (9) numerically. This process is the SDM's local analysis procedure. A standard FDM solver is used in this study.

First, we divide $L$ into a grid, as shown in Fig. 3(a). The grid interval is $h$, and $u_{i,j}$ and $f_{i,j}$ denote the variable and the source, respectively, at the position
$$\mathbf{x}_{i,j} = (ih - 0.5h, jh - 0.5h)^T. \quad (20)$$
These parameters are given by
$$u_{i,j} = u(\mathbf{x}_{i,j}), \quad (21)$$
$$f_{i,j} = f(\mathbf{x}_{i,j}). \quad (22)$$

Let us now consider a case where all $f_{i,j}$ are known for all grid points in $L$.

Separately from $u_{i,j}$, we define the temperature at the $k$th grid point on $\partial L$ as $u_k = u(s_k)$ (see Fig. 3(b)). As stated in Subsection 2.2.2.1 and Eq. (12), $u_k$ is expressed as a product of $\mathbf{N}^{\partial L}$ and $\mathbf{u}^L = (u_{\mathrm{CP1}} \cdots u_{\mathrm{CP4}})^T$:
$$u_k = u(s_k) = \mathbf{N}^{\partial L}(s_k)\mathbf{u}^L. \quad (23)$$

The process of construction of the FDM model of $L$ generates approximations for $b$ and $c$ (which are referred to as $b_{\mathrm{FDM}}$ and $c_{\mathrm{FDM}}$, respectively). $b_{\mathrm{ERROR}}$ and $c_{\mathrm{ERROR}}$ denote the differences between the true functions $b_{\mathrm{TRUE}}$ and $c_{\mathrm{TRUE}}$ and the approximations $b_{\mathrm{ERROR}}$ and $c_{\mathrm{ERROR}}$, respectively. We then obtain
$$\begin{aligned} b &= b_{\mathrm{FDM}} = b_{\mathrm{TRUE}} + b_{\mathrm{ERROR}} \\ c &= c_{\mathrm{FDM}} = c_{\mathrm{TRUE}} + c_{\mathrm{ERROR}}. \end{aligned} \quad (24)$$
The effects of $b_{\mathrm{ERROR}}, c_{\mathrm{ERROR}}$ on the accuracy of $u_{\mathrm{CP5}}$ are explained in Subsection 2.2.2.3.

### 2.2.2.3 Summary of errors in the SDM

This subsection deals with how the error factors of the SDM that were stated in Subsections 2.2.2.1 and 2.2.2.2 affect the precision of $u_{\mathrm{CP5}}$. When Eq. (24) is substituted into Eq. (18), it follows that

$$\begin{aligned} u_{\mathrm{CP5}} &= \int_{\partial L}(b_{\mathrm{TRUE}} + b_{\mathrm{ERROR}})uds \\ &\quad + \int_L (c_{\mathrm{TRUE}} + c_{\mathrm{ERROR}})f d\mathbf{x} \\ &= \int_{\partial L}(b_{\mathrm{TRUE}} + b_{\mathrm{ERROR}})(\mathbf{N}^{\partial L}_{\mathrm{TRUE}} + \mathbf{N}^{\partial L}_{\mathrm{ERROR}})ds\mathbf{u}^L \\ &\quad + \int_L (c_{\mathrm{TRUE}} + c_{\mathrm{ERROR}})f d\mathbf{x} \\ &= (\mathbf{a}^L_{\mathrm{TRUE}} + \mathbf{a}^L_{\mathrm{ERROR}})\mathbf{u}^L + F^L_{\mathrm{TRUE}} + F^L_{\mathrm{ERROR}} \\ &= \mathbf{a}^L \mathbf{u}^L + F^L(\mathbf{x}_{\mathrm{CP5}}), \end{aligned} \quad (25)$$

where
$$\begin{aligned} \mathbf{a}^L_{\mathrm{TRUE}} &= \int_{\partial L} b_{\mathrm{TRUE}}\mathbf{N}^{\partial L}_{\mathrm{TRUE}}ds \\ \mathbf{a}^L_{\mathrm{ERROR}} &= \int_{\partial L} (b_{\mathrm{TRUE}}\mathbf{N}^{\partial L}_{\mathrm{ERROR}} + b_{\mathrm{ERROR}}\mathbf{N}^{\partial L}_{\mathrm{TRUE}} \\ &\quad + b_{\mathrm{ERROR}}\mathbf{N}^{\partial L}_{\mathrm{ERROR}})ds \\ F^L_{\mathrm{TRUE}} &= \int_L c_{\mathrm{TRUE}} f d\mathbf{x} \\ F^L_{\mathrm{ERROR}} &= \int_{\Omega_L} c_{\mathrm{ERROR}} f d\mathbf{x}, \end{aligned} \quad (26)$$

and
$$\begin{aligned} \mathbf{a}^L &= \mathbf{a}^L_{\mathrm{TRUE}} + \mathbf{a}^L_{\mathrm{ERROR}} \\ F^L(\mathbf{x}_{\mathrm{CP5}}) &= F^L_{\mathrm{TRUE}} + F^L_{\mathrm{ERROR}}. \end{aligned} \quad (27)$$

If we can prepare a completely accurate $\mathbf{N}^{\partial L}$ (which is the interpolation matrix for $\partial L$), the following conditions are then satisfied:
· $\mathbf{N}^{\partial L}_{\mathrm{ERROR}} = 0$
· The error of the influence coefficient matrix $\mathbf{a}^L$, denoted by $\mathbf{a}^L_{\mathrm{ERROR}}$, is expressed in the following form:
$$\mathbf{a}^L_{\mathrm{ERROR}} = \int_{\partial L} b_{\mathrm{ERROR}} \mathbf{N}^{\partial L}_{\mathrm{TRUE}} ds. \quad (28)$$
$\mathbf{a}^L_{\mathrm{ERROR}}$ therefore only includes an error related to $b_{\mathrm{ERROR}}$.
· When both $\mathbf{N}^{\partial L}$ and $\mathbf{u}^L$ are correct, the precise variable profile on $\partial L$ (where $u_k = u(s_k)$) can be derived. If $\mathbf{N}^{\partial L}$ is accurate for all local domains, then the variable values of the grid points on the boundary of every local domain are exactly equivalent to those of its neighboring local domain. In this case, the SDM solution corresponds exactly to that of the direct FDM model, which has a grid that is as fine as that of $L$ throughout the global field. In this case, when analyzing $L$, which has a grid interval of $h$, using the second-order finite-difference discretization, the order of the error of the SDM is $O(h^2)$. The SDM's error is at the same level as that obtained from the direct FDM model, which includes the entire simulated field and where the grid interval is $h$ throughout the field. Basically, this indicates that the accuracy of the SDM does not exceed that of the direct FDM.

## 3 Computational implementation of the SDM using FDM in the local analysis

In the numerical experiments in Sections 4 and 5, we simulate the 2D linear stationary temperature fields without a heat source. The source term on the right side of the governing equation (Eq. (1)) is zero (i.e., $f(\mathbf{x}) = 0$) in this case. This means that

$$\frac{\partial}{\partial x}\left(k\frac{\partial u}{\partial x}\right) + \frac{\partial}{\partial y}\left(k\frac{\partial u}{\partial y}\right) = 0 \text{ in } G. \quad (29)$$

We describe the computational implementation of the SDM for the governing equation above in this section.

## 3.1 Local analysis

The objective of the local analysis process is to generate the interpolation matrix $\mathbf{N}^L$ and the influence coefficient matrix $\mathbf{a}^L$. The local analysis is conducted using the FDM (Subsection 2.2) in the following manner.

$L$ in Fig. 2(a) shows CP 5 at its center and reference CPs 1–4 are arranged around CP 5. The vector (including the temperature values) at CPs 1–4 is given by $\mathbf{u}^L = (u_{CP1} \cdots u_{CP4})^T$.

We spatially discretize $L$ using a grid interval of $h$ to construct the FDM model. $u_{i,j}$ denotes the temperature at position $\mathbf{x}_{i,j} = (ih - 0.5h, jh - 0.5h)^T$. Because no heat sources are located within the domain ($f(\mathbf{x}) = 0$), $u_{i,j}$ in Eq. (9) can be expressed discretely using the form

$$u_{i,j} = u(\mathbf{x}_{i,j}) = \mathbf{N}^L(\mathbf{x}_{i,j})\mathbf{u}^L$$
$$\mathbf{N}^L(\mathbf{x}_{i,j}) = \int_{\partial L} b(\mathbf{x}_{i,j}, s)\mathbf{N}^{\partial L}(s)ds \quad (30)$$
$$= \sum_k w(s_k) b(\mathbf{x}_{i,j}, s_k)\mathbf{N}^{\partial L}(s_k).$$

where $w(s_k)$ is the weighting factor at the position $s_k$.

Based on substitution of $\mathbf{x}_{i,j} = \mathbf{x}_{CP5}$, we can then obtain $\mathbf{a}^L$.

$$u_{CP5} = u(\mathbf{x}_{CP5}) = \mathbf{N}^L(\mathbf{x}_{CP5})\mathbf{u}^L \quad (31)$$
$$\mathbf{a}^L = \mathbf{N}^L(\mathbf{x}_{CP5}),$$

Next, we must introduce two detailed ways to construct $\mathbf{N}^L$ and $\mathbf{a}^L$. The first is to perform a local analysis without the oversampling technique [9,10] (which will be described in Subsection 3.1.1) and the other is to perform a local analysis using the oversampling scheme (which will be described in Subsection 3.1.2). The latter method generates a more accurate $\mathbf{N}^L$ in many cases.

### 3.1.1 Local analysis without oversampling

(1) Begin by constructing fine-grained finite-difference local domains using a grid interval of $h$.

(2) Then, determine the boundary condition of $L$, i.e., $u(s)$ on $\partial L$.

The position coordinates of the grid points on the closed curve $\partial L$ are denoted by $s_1, s_2 \cdots$. The temperatures for the boundary grid points are thus denoted by $u(s_1), u(s_2) \cdots$, as indicated in Fig. 3(b).

Similar to Eq. (12) in Subsection 2.2.2.1, $u(s_i)$ is computed by interpolation of the temperature for the CPs on $\partial L$, i.e.,

$$\mathbf{u}^L = (u_{CP1} \cdots u_{CP4})^T.$$

$$u(s_i) = \mathbf{N}^{\partial L}(s_i)\mathbf{u}^L \text{ on } \partial L. \quad (32)$$

The simplest form of interpolation is the linear form. For example, the temperatures for the grid points that are located between CPs 1 and 2 are expressed using a linear interpolation of $u_{CP1}$ and $u_{CP2}$.

When the temperature distribution in $L$ is sufficiently gradual and when $L$ is homogeneous (i.e., the conductivity is uniform throughout $L$), the linear interpolation may generate sufficiently accurate $u(s_i)$ values. Otherwise, the linear interpolation leads to an error. For example, if $L$ is heterogeneous, $u(s_i)$ will have a complex distribution based on the nonuniform conductivity profile of $L$. In this case, $u(s_i)$ would not have a linear distribution.

As shown in Eq. (19) in Subsection 2.2.2.1, the error in $u(s_i)$ is a result of $\mathbf{N}^{\partial L}_{ERROR}$ (which is the error of $\mathbf{N}^{\partial L}$). Specifically, the error for $u(s_i)$ is equivalent to the second term, $\mathbf{N}^{\partial L}_{ERROR}(s_i)\mathbf{u}^L$, in the equation below.

$$u(s_i) = \mathbf{N}^{\partial L}(s_i)\mathbf{u}^L$$
$$= \left(\mathbf{N}^{\partial L}_{TRUE}(s_i) + \mathbf{N}^{\partial L}_{ERROR}(s_i)\right)\mathbf{u}^L \quad (33)$$
$$= \mathbf{N}^{\partial L}_{TRUE}(s_i)\mathbf{u}^L + \mathbf{N}^{\partial L}_{ERROR}(s_i)\mathbf{u}^L.$$

(3) Formulate the finite difference equations for $L$ to generate $\mathbf{N}^L$;

By formulating the finite difference equations for $L$ by imposing $u(s_i)$ from Eq. (32) as the boundary conditions on $L$, we obtain an approximation of $\mathbf{N}^L$ in Eq. (30).

(4) Finally, derive the influence coefficient matrix, $\mathbf{a}^L = \mathbf{N}^L(\mathbf{x}_{CP5})$, from Eq. (17)

### 3.1.2 Local analysis with oversampling

(1) Begin by constructing FDM models of the local domains, including the oversampled domain $L^+$.

To reduce the error in $\mathbf{N}^{\partial L}$ (i.e., the $\mathbf{N}^{\partial L}_{ERROR}$ shown in Subsections 2.2.2.1 and 2.2.2.3), we use an oversampling technique [9,10] in the local analysis. As shown in Fig. 4, we extract a local domain, $L^+$, that is larger than $L$. $L^+$ is composed of $L$ and its surrounding domain (which is the oversampled domain). We construct the FDM model of $L^+$.

(2) Then, determine the boundary condition of $L^+$, i.e., $u(s_i^+)$ on $\partial L^+$.

The outermost boundary of $L^+$ is $\partial L^+$. To determine the boundary condition, $u(s_i^+)$ on $\partial L^+$, we follow the procedures below.

We define the grid points at the four corners of $L^+$ as the boundary points (shown as triangles in Fig. 4(b)). The temperatures at these points are denoted by

$$\mathbf{u}^{L+} = (u_1^{L+} \cdots u_4^{L+})^T. \quad (34)$$

In this manuscript, we linearly interpolate $\mathbf{u}^{L+}$ to generate $u(s_i^+)$ on $\partial L^+$.

$$u(s_i^+) = \mathbf{N}^{\partial L+}(s_i^+)\mathbf{u}^{L+} \text{ on } \partial L^+. \quad (35)$$

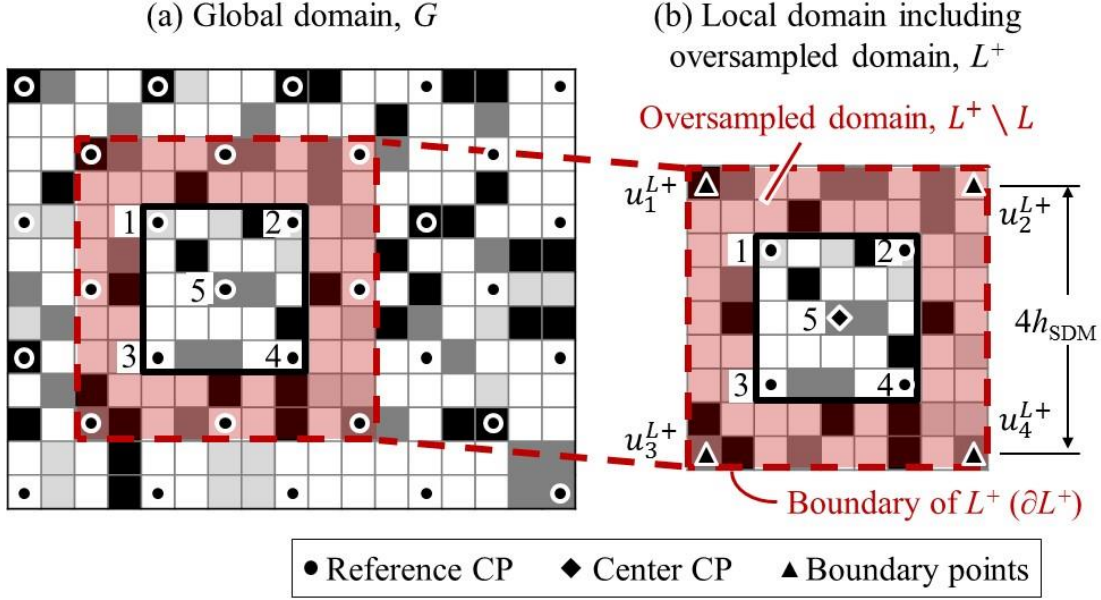

Fig. 4. (a) Nonperiodic global domain, and (b) local domain including an oversampled domain, $L^+$.

where $u(s_i^+)$ is the temperature for the $i$th grid point on $\partial L^+$ and $\mathbf{N}^{\partial L^+}$ is the interpolation function matrix for $\partial L^+$.

(3) Formulate the finite difference equations for $L^+$ to construct $\mathbf{N}^L$.

We formulate the finite difference equations for $L^+$ by imposing $u(s_i^+)$ (from Eq. (35)) as the boundary conditions on $L^+$ and then solve them. We then obtain
$$u_{i,j} = u(\mathbf{x}_{i,j}) = \mathbf{N}^{L^+}(\mathbf{x}_{i,j})\mathbf{u}^{L^+} \text{ in } L^+, \quad (36)$$
where $\mathbf{N}^{L^+}$ is the interpolation matrix for $L^+$.

$u_{\text{CP}i}$ is then computed by substituting $\mathbf{x}_{i,j} = \mathbf{x}_{\text{CP}i}$ into the above equation to give
$$u_{\text{CP}i} = u(\mathbf{x}_{\text{CP}i}) = \mathbf{N}^{L^+}(\mathbf{x}_{\text{CP}i})\mathbf{u}^{L^+}. \quad (37)$$
From the equation above, we obtain the following two relations.
$$u_{\text{CP}5} = u(\mathbf{x}_{\text{CP}5}) = \mathbf{N}^{L^+}(\mathbf{x}_{\text{CP}5})\mathbf{u}^{L^+}, \quad (38)$$
$$\mathbf{u}^L = \mathbf{d}^{L^+}\mathbf{u}^{L^+}, \quad (39)$$
where
$$\mathbf{d}^{L^+} = (\mathbf{N}^{L^+}(\mathbf{x}_{\text{CP}1})^T \cdots \mathbf{N}^{L^+}(\mathbf{x}_{\text{CP}4})^T)^T. \quad (40)$$
If $\mathbf{d}^{L^+}$ is a regular matrix, we can then calculate its inverse matrix, $(\mathbf{d}^{L^+})^{-1}$; i.e.,
$$\mathbf{u}^{L^+} = (\mathbf{d}^{L^+})^{-1}\mathbf{u}^L. \quad (41)$$
By substituting this into Eq. (36), we then obtain
$$u_{i,j} = u(\mathbf{x}_{i,j}) = \mathbf{N}^{L^+}(\mathbf{x}_{i,j})\mathbf{u}^{L^+}$$
$$= \mathbf{N}^{L^+}(\mathbf{x}_{i,j})(\mathbf{d}^{L^+})^{-1}\mathbf{u}^L = \mathbf{N}^L(\mathbf{x}_{i,j})\mathbf{u}^L \text{ in } L^+, \quad (42)$$
where
$$\mathbf{N}^L(\mathbf{x}_{i,j}) = \mathbf{N}^{L^+}(\mathbf{x}_{i,j})(\mathbf{d}^{L^+})^{-1}, \quad (43)$$
is the interpolation matrix that estimates the temperature at an arbitrary grid point denoted by $(i,j)$ in $L^+$ with reference to $\mathbf{u}^L$.

(4) Finally, determine the influence coefficient matrix $\mathbf{a}^L$ by substituting $\mathbf{x} = \mathbf{x}_{\text{CP}5}$ into Eq. (43).
$$\mathbf{a}^L = \mathbf{N}^L(\mathbf{x}_{\text{CP}5}) = \mathbf{N}^{L^+}(\mathbf{x}_{\text{CP}5})(\mathbf{d}^{L^+})^{-1}. \quad (44)$$
Here, we consider the reason why the oversampling technique should be used. The important points are written below.

・$L$ does not include the oversampled domain and has the boundary $\partial L$.

・$L^+$ does include the oversampled domain and has the boundary $\partial L^+$.

When the oversampling scheme is used, we impose the Dirichlet boundary condition, i.e., $u(s_i^+) = \mathbf{N}^{\partial L^+}(s_i^+)\mathbf{u}^{L^+}$, on $\partial L^+$ and then conduct the FDM analysis of $L^+$.

By substituting the position of the $i$th grid point on $\partial L$ (i.e., NOT on $\partial L^+$) into Eq. (42), we then obtain
$$u(s_i) = \mathbf{N}^{L^+}(s_i)(\mathbf{d}^{L^+})^{-1}\mathbf{u}^L \text{ on } \partial L. \quad (45)$$
We now focus on the difference in $u(s_i)$ due to the use of the oversampling technique.

・Without the oversampling technique: $u(s_i) = \mathbf{N}^{\partial L}(s_i)\mathbf{u}^L$ on $\partial L$;

・With the oversampling technique:
$u(s_i) = \mathbf{N}^{L^+}(s_i)(\mathbf{d}^{L^+})^{-1}\mathbf{u}^L$ on $\partial L$.

Based on the above, the effects of the oversampling are described as follows. The temperature profile in $L$ in the case where we apply the boundary condition (i.e., where $u(s_i^+) = \mathbf{N}^{\partial L^+}(s_i^+)\mathbf{u}^{L^+}$) on $\partial L^+$ is equivalent to that in the case where $u(s_i) = \mathbf{N}^{L^+}(s_i)(\mathbf{d}^{L^+})^{-1}\mathbf{u}^L$ is imposed on $\partial L$.

$u(s_i^+)$ on $\partial L^+$ is simply a temperature profile that has been linearly interpolated from the temperatures for boundary points 1–4 on $\partial L^+$. However, this attempts to make $u(s_i)$ valid based on the nonuniform thermal conductivity distributions that occur both inside and outside $L$, which is a benefit of the oversampling technique.

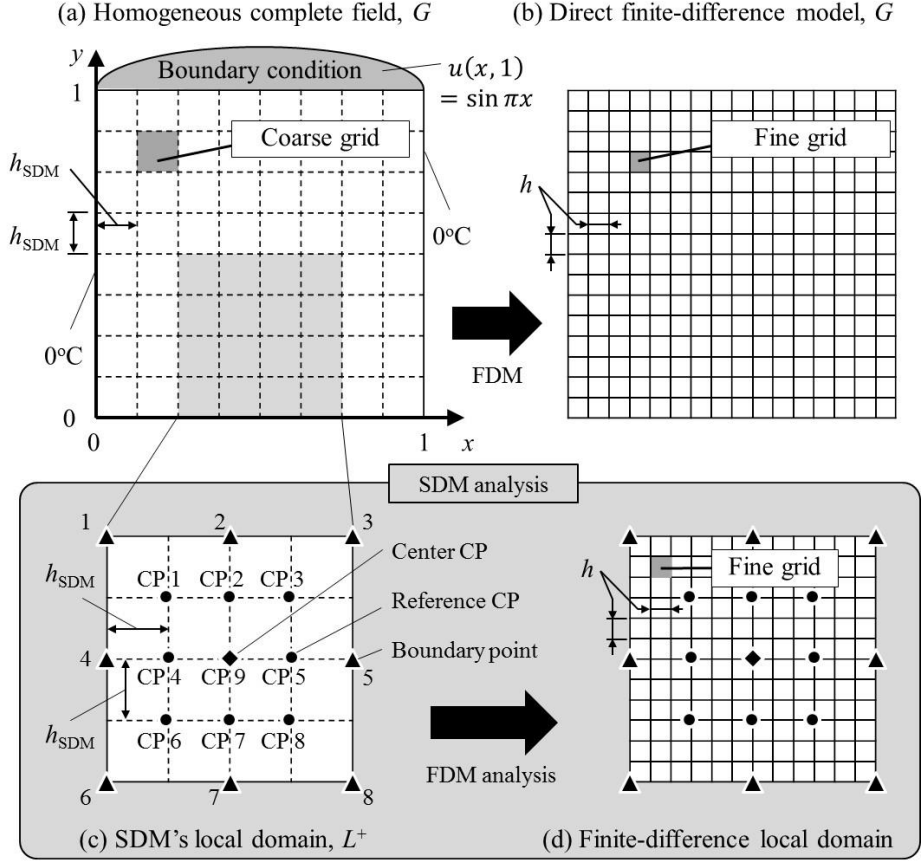

Fig. 5. (a) Homogeneous and isotropic global domain (two-dimensional stationary temperature field) for the convergence investigation; (b) direct FDM model; (c) local domain of the SDM, $L^+$; and (d) finite-difference model of $L^+$.

### 3.2 Global analysis

We move on to the global analysis after construction of the influence coefficient matrix $\mathbf{a}^L$ and the interpolation function matrix $\mathbf{N}^L$ for all local domains over the entire analytical field $G$, using the analysis steps below.

(1) Arrange the CPs in $G$.
We arrange $n$ CPs (designated CPs 1–$n$) in $G$, where $n$ is an arbitrary number.
$$\mathbf{u}^G = (u^G_{CP1} \cdots u^G_{CPn})^T \quad \text{in } G, \tag{46}$$
where $u^G_{CPi}$ is the temperature value of CP $i$.

(2) Construct relational expressions for the CPs in $G$.
We formulate relational expressions for the CPs (similar to Eq. (31)) for all local domains in $G$. We then obtain
$$\mathbf{u}^G = \mathbf{a}^G \mathbf{u}^G \quad \text{in } G, \tag{47}$$
where $\mathbf{a}^G$ is the global influence coefficient matrix. $\mathbf{a}^G$ is then established by assembling the entries from all local influence coefficient matrices ($\mathbf{a}^L$). $\mathbf{a}^G$ is a band matrix with a bandwidth that is approximately the same as the number of reference CPs. In $L$, as shown in Sections 2 and 3, the temperature at the center CP is computed with reference to the four reference CPs. The bandwidth in this case is thus approximately four.

(3) Solve Eq. (47) to obtain $\mathbf{u}^G$.
We cannot solve Eq. (47) in its current form because Eq. (47) is an identity function of $\mathbf{u}^G$. By applying the global boundary conditions, we can then solve Eq. (47) and subsequently compute $\mathbf{u}^G$.

(4) Interpolate the local temperature distributions.
A detailed temperature distribution is obtained for each local domain by interpolation of the global temperature solution ($\mathbf{u}^G$) using $\mathbf{N}^L$. The nonuniform conductivity profile is taken into account during this interpolation process. The temperature profile throughout $G$ can then be generated by connecting all the local profiles.

## 4 Convergence study: Homogeneous temperature field

All numerical experiments in this paper were performed using a workstation with the following specification:
- Intel® Core i7-3930K central processing unit (CPU; 3.20 GHz, six cores, 612 threads);
- 64 GB of random access memory (RAM).

We use MATLAB R2016a (MathWorks, Inc.) and set the format for the numerical values to a 15-digit scaled fixed point format.

## 4.1 Problem statement

- The target field (global domain $G$) is a 2D square, as shown in Fig. 5(a), with side length of 1.0 (dimensionless length);
- $G$ is a linear steady-state temperature field;
- The thermal conductivity profile throughout $G$ is both uniform and isotropic. The governing equation is thus
$$\frac{\partial^2 u}{\partial x^2} + \frac{\partial^2 u}{\partial y^2} = 0, \quad (48)$$
where $u = u(\mathbf{x})$ is the temperature at $\mathbf{x} = (x, y)^T$;
- The Dirichlet boundary conditions (which are fixed-temperature boundary conditions) that are shown in Fig. 5(a) are imposed on $G$. The temperature profile on the top side of $G$, denoted by $u(x, 1)$, is the half-wave length of a sine wave, i.e.,
$$u(x, 1) = sin(\pi x).$$
The temperature on the other three sides is fixed at 0°C;
- This problem can be solved exactly and the true temperature is given as follows:
$$u_{\text{TRUE}}(\mathbf{x}) = u_{\text{TRUE}}(x, y) = sin(\pi x) \frac{exp(\pi y) - exp(-\pi y)}{exp(\pi) - exp(-\pi)}. \quad (49)$$
- This problem is then solved using the SDM, the direct FDM, and the DDM. We then use the temperature difference when compared with the true solution as an index of accuracy and can thus compare the convergence properties of the three methods.

## 4.2 Direct FDM

$G$ is divided into an $n_{\text{grid}}^G \times n_{\text{grid}}^G$ grid, as shown in Fig. 5(b). The grid interval $h$ is given by:
$$h = \frac{1}{n_{\text{grid}}^G - 1}. \quad (50)$$

The total number of grid points in the direct FDM model of $G$ is given by:
$$n_{\text{FDM}}^G = \left(n_{\text{grid}}^G - 1\right)^2. \quad (51)$$

We then test the following four FDM models with the different values of $h$ below.
$$h = \frac{1}{8}, \frac{1}{16}, \frac{1}{32}, \frac{1}{64}. \quad (52)$$

We first derive the finite difference equation for the grid point $(i, j)$. By discretizing Eq. (48) using a second-order central difference, we then obtain:
$$u_{i,j} = \frac{u_{i-1,j} + u_{i+1,j} + u_{i,j-1} + u_{i,j+1}}{4} + O(h^2). \quad (53)$$
The order of the error of $u_{i,j}$ is $O(h^2)$.

In the direct FDM analysis, we must solve a linear algebraic equation with a square matrix of order $n_{\text{FDM}}^G$ and a bandwidth of 4.

## 4.3 SDM
### 4.3.1 Global analysis

The global analysis is conducted in accordance with the protocols that were described in Subsection 3.2. We arrange the CPs at equal intervals of $h_{\text{SDM}}$ on $G$, as shown in Fig. 5(a). We then prepare three SDM models using the different $h_{\text{SDM}}$ values below:
$$h_{\text{SDM}} = \frac{1}{8}, \frac{1}{16}, \frac{1}{32}. \quad (54)$$
The number of CPs on one side of $G$ is denoted by $n_{\text{CP}}^G$, which is determined as follows:
$$n_{\text{CP}}^G = \frac{h}{h_{\text{SDM}}}\left(n_{\text{grid}}^G - 1\right) + 1 \approx \frac{h}{h_{\text{SDM}}} n_{\text{grid}}^G. \quad (55)$$
Because $G$ is a square, the total number of CPs in $G$, denoted by $n_{\text{SDM}}^G$, is:
$$n_{\text{SDM}}^G = (n_{\text{CP}}^G)^2 = \left(\frac{h}{h_{\text{SDM}}}\left(n_{\text{grid}}^G - 1\right) + 1\right)^2 \approx \frac{h^2}{h_{\text{SDM}}^2} n_{\text{FDM}}^G. \quad (56)$$
When $n_{\text{grid}}^G$ is sufficiently large, $n_{\text{CP}}^G$ is then $h/h_{\text{SDM}}$ times larger than $n_{\text{grid}}^G$. Therefore, $n_{\text{SDM}}^G$ can be suppressed to a value that is $(h/h_{\text{SDM}})^2$ times as many as that of $n_{\text{FDM}}^G$.

The global influence coefficient matrix $\mathbf{a}^G$ is established by assembling all the influence coefficient matrices, $\mathbf{a}^L$. We then construct relational expressions for all the CPs in $G$, where $\mathbf{u}^G = \mathbf{a}^G \mathbf{u}^G$, and then solve them using the global boundary conditions. The order and the bandwidth of $\mathbf{a}^G$ are $n_{\text{SDM}}^G$ and 8, respectively.

### 4.3.2 Local analysis

To construct $\mathbf{a}^L$ and $\mathbf{N}^L$, the local analysis is performed using the oversampling scheme according to the procedures that were introduced in Subsection 3.1.2. However, unlike the scenario in Subsection 3.1.2, we arrange for eight reference CPs to be in each local domain $L^+$. Therefore, the temperature at the center at CP 9 ($u_{\text{CP9}}$) is expressed as a product of $\mathbf{a}^L$ and the temperatures at the eight reference CPs, given by $\mathbf{u}^L = (u_{\text{CP1}} \cdots u_{\text{CP8}})^T$, as shown in Fig. 5(c).
$$u_{\text{CP9}} = \mathbf{a}^L \mathbf{u}^L. \quad (57)$$
In the FDM analysis of each of the local domains, we use the same second-order finite-difference discretization (see Subsection 4.2 and Eq. (53)) that was used in the direct FDM. The grid interval $h$ in this case is
$$h = \frac{1}{64}. \quad (58)$$
$L^+$ is a square with a side length of $4h_{\text{SDM}}$ (i.e., it is four times larger than the interval between the CPs, $h_{\text{SDM}}$). The total number of grid points in $L^+$ is denoted by $n_{\text{SDM}}^{L^+}$, where:
$$n_{\text{SDM}}^{L^+} = \left(4\frac{h_{\text{SDM}}}{h} + 1\right)^2. \quad (59)$$

From the above, in the FDM analysis of $L^+$, we must solve a linear algebraic equation with a square matrix of order $n_{\text{SDM}}^{L^+}$ and a bandwidth of 4.

$G$ in this section is homogeneous and thus there is no difference in thermal conductivity between the local domains. In this case, an influence coefficient matrix that has been derived for one local domain is applicable to all other local domains. Therefore, we only need to perform the local analysis once.

Based on the above, we must solve a linear algebraic equation with a square matrix of order $n_{\text{SDM}}^{L^+}$ only once in the local analysis.

### 4.4 DDM

The detailed computational protocols for the DDM are given in Subsection 5.4 and are thus excluded here.

In a similar manner to the SDM, $G$ is divided into smaller local domains, i.e., square domains with side lengths of $h_{\text{DDM}}$, in $L_{\text{DDM}}$. We then prepare three DDM models with different values of $h_{\text{DDM}}$:
$$h_{\text{DDM}} = \frac{1}{8}, \frac{1}{16}, \frac{1}{32}. \tag{60}$$
In the same manner as the SDM, we use the FDM solver to perform the DDM's local analysis. The grid interval used for the finite difference local models is $h = 1/64$, which is the same as that used in the finest model of the four direct FDM models that were described in Subsection 4.2.

### 4.5 Results

The exact solution to this problem is given in Eq. (49). The temperature values that were obtained from the direct FDM models, the SDM models, and the DDM models are denoted by $u_{\text{FDM}}$, $u_{\text{SDM}}$, and $u_{\text{DDM}}$, respectively. We define their root mean square errors (RMSEs) using the following forms and show these RMSEs in Fig. 6.

$$\begin{aligned}
RMSE_{\text{FDM}} &= \sqrt{\frac{1}{n_{\text{FDM}}^G} \sum_{i,j} \left(u_{\text{FDM}}(\mathbf{x}_{i,j}) - u_{\text{TRUE}}(\mathbf{x}_{i,j})\right)^2} \\
RMSE_{\text{SDM}} &= \sqrt{\frac{1}{n_{\text{SDM}}^G} \sum_{i,j} \left(u_{\text{SDM}}(\mathbf{x}_{i,j}) - u_{\text{TRUE}}(\mathbf{x}_{i,j})\right)^2} \\
RMSE_{\text{DDM}} &= \sqrt{\frac{1}{n_{\text{DDM}}^G} \sum_{i,j} \left(u_{\text{DDM}}(\mathbf{x}_{i,j}) - u_{\text{TRUE}}(\mathbf{x}_{i,j})\right)^2}.
\end{aligned} \tag{61}$$

The vertical axis indicates the logarithmic RMSEs, while the horizontal axis shows the distances between two adjacent discretization points for each method (i.e., direct FDM: $h$; SDM: $h_{\text{SDM}}$; DDM: $h_{\text{DDM}}$).

#### 4.5.1 FDM

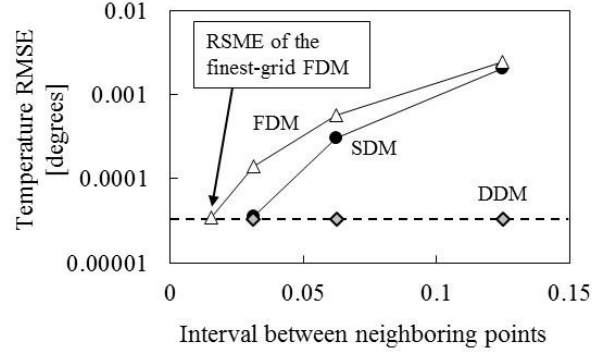

Fig. 6. Results of the convergence investigation, including temperature RMSEs of the direct FDM models, the SDM models, and the DDM models.

The $RMSE_{\text{FDM}}$ is proportional to the squared grid interval. This is a reasonable result because a second-order finite difference discretization has been used.

#### 4.5.2 SDM

The grid interval used for the finite difference local domains is $h = 1/64$, which is the same as that used for the finest model of the four direct FDM models that were described in Subsection 4.2. The SDM solution only converges to the solution of the finest direct FDM (the dashed line in Fig. 6) when $\mathbf{a}_{\text{ERROR}}^L$ in Eq. (27) is a zero matrix for all local domains.

$RMSE_{\text{SDM}}$ appears to be proportional to the squared interval between the CPs. When the number of divided regions used is the same, $RMSE_{\text{SDM}}$ is always smaller than $RMSE_{\text{FDM}}$. In the case where $h_{\text{SDM}} = 1/32$, the SDM solution converges to the finest FDM solution (the dashed line in Fig. 6).

#### 4.5.3 DDM

Like the SDM, the grid interval used for the DDM's local domains is the same as that used for the finest direct FDM ($h = 1/64$). $RMSE_{\text{DDM}}$ is almost the same as the $RMSE_{\text{FDM}}$ of the finest FDM (the dashed line in Fig. 6), regardless of the number of divided regions used. Unlike the SDM, the grid points on the boundaries of the local domains in the DDM are not coarse-grained (reduced). This guarantees that the DDM can solve problems with the same level of accuracy as the direct FDM.

## 5 Investigation of practical feasibility: Nonuniform and nonperiodic temperature field

As stated earlier, all numerical experiments in this work were performed on a workstation with the following specification:
- Intel® Core i7-3930K CPU (3.20 GHz, six cores, 612

- threads);
- 64 GB of RAM.

We used MATLAB R2016a and set the format for the numerical values to 15-digit scaled fixed point. We do not use parallel computation methods for the local analyses.

## 5.1 Problem statement

- The target field (global domain $G$) is the 2D square shown in Fig. 7, which has a side length of 1.0 (dimensionless length);
- $G$ has a nonuniform and nonperiodic thermal conductivity distribution; the governing equation for $G$ is

$$0 = \frac{\partial}{\partial x}\left(k\frac{\partial u}{\partial x}\right) + \frac{\partial}{\partial y}\left(k\frac{\partial u}{\partial y}\right) \text{ in } G, \quad (62)$$

where $u = u(\mathbf{x})$ and $k = k(\mathbf{x})$ are the temperature and the conductivity at $\mathbf{x} = (x,y)^T$, respectively;

- $G$ consists of $(n_{\text{grid}}^G - 1) \times (n_{\text{grid}}^G - 1)$ small squares with the same side length of $h = 1/(n_{\text{grid}}^G - 1)$;
- The small squares have different conductivities. The conductivity is constant in each small square and is isotropic; the conductivity in square $(i,j)$ (i.e., the square at the $i$th row in the $j$th column) is denoted by $k_{i,j}$.

$$k(\mathbf{x}) = k(x,y) = k_{i,j} \text{ for } x \in (ih-h, ih), y \in (jh-h, jh), \quad (63)$$

$k_{i,j}$ values are uniformly distributed between 0.01–1.00; the difference in conductivity between the small squares is thus at most 100 times.

$$0.01 \leq k_{i,j} \leq 1.00. \quad (64)$$

Note that $k_{i,j}$ represents the dimensionless conductivity;

- The global boundary conditions are shown in Fig. 7. The temperatures at the four corners of $G$ are set, and the temperature profiles on the four sides are linearly interpolated based on the four corners;

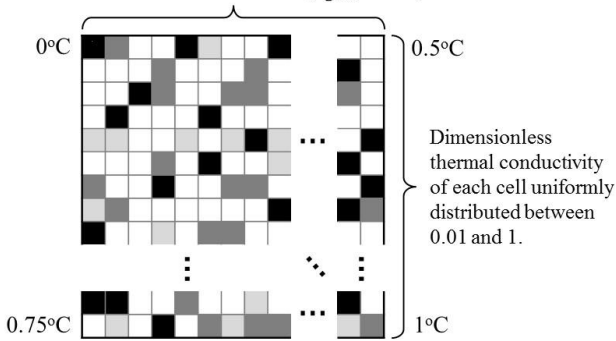

Fig. 7. Example of nonperiodic and heterogeneous global domain problem with two-dimensional stationary temperature field.

- This problem cannot be solved exactly; we therefore compare the direct FDM with the SDM and the DDM in terms of both computational precision and time requirements.

## 5.2 Direct FDM

$G$ is divided into a $n_{\text{grid}}^G \times n_{\text{grid}}^G$ grid with an interval of $h = 1/(n_{\text{grid}}^G - 1)$. Therefore, each small square has a single grid point at its center. The total number of grid points contained in the direct FDM model of $G$ is $n_{\text{FDM}}^G = (n_{\text{grid}}^G - 1)^2$.

To derive the finite difference equation, we arrange the virtual grid points at the midpoint of the grid points, which is simply called the middle grid point. For example, the middle grid point $(i + 1/2, j)$ is located between grid points $(i,j)$ and $(i+1,j)$. When the second-order partial derivative of $u$ with respect to $x$ is expressed using the middle grid points, it follows that

$$\frac{\partial}{\partial x}\left(k\frac{\partial u}{\partial x}\right)_{i,j} = \frac{\left(k\frac{\partial u}{\partial x}\right)_{i+1/2,j} - \left(k\frac{\partial u}{\partial x}\right)_{i-1/2,j}}{h} + O(h^2). \quad (65)$$

where the two first-order partial derivatives of $u$ on the right side are given as follows:

$$\begin{aligned}\left(k\frac{\partial u}{\partial x}\right)_{i+1/2,j} &= k(\mathbf{x}_{i+1/2,j})\left(\frac{\partial u}{\partial x}\right)_{i+1/2,j} \\ &= k_{i+1/2,j}\frac{u_{i+1,j} - u_{i,j}}{h} + O(h^2) \\ \left(k\frac{\partial u}{\partial x}\right)_{i-1/2,j} &= k(\mathbf{x}_{i-1/2,j})\left(\frac{\partial u}{\partial x}\right)_{i-1/2,j} \\ &= k_{i-1/2,j}\frac{u_{i,j} - u_{i-1,j}}{h} + O(h^2),\end{aligned} \quad (66)$$

where

$$k_{i+1/2,j} = k(\mathbf{x}_{i+1/2,j}), k_{i-1/2,j} = k(\mathbf{x}_{i-1/2,j}), \quad (67)$$

are the equivalent thermal conductivities that are given in Eqs. (74) and (75) below. We then substitute Eq. (66) into Eq. (65) to obtain

$$\frac{\partial}{\partial x}\left(k\frac{\partial u}{\partial x}\right)_{i,j} = \frac{k_{i+1/2,j}(u_{i+1,j}-u_{i,j}) - k_{i-1/2,j}(u_{i,j}-u_{i-1,j})}{h^2} + O(h). \quad (68)$$

In a similar manner to the derivation with respect to $x$, we can also derive the second-order partial derivative of $u$ with respect to $y$:

$$\frac{\partial}{\partial y}\left(k\frac{\partial u}{\partial y}\right)_{i,j} = \frac{k_{i,j+1/2}(u_{i,j+1}-u_{i,j}) - k_{i,j-1/2}(u_{i,j}-u_{i,j-1})}{h^2} + O(h). \quad (69)$$

From the above, we obtain

$$\begin{aligned}0 &= \frac{\partial}{\partial x}\left(k\frac{\partial u}{\partial x}\right)_{i,j} + \frac{\partial}{\partial y}\left(k\frac{\partial u}{\partial y}\right)_{i,j} = \\ &\frac{1}{h^2}\Big(k_{i+1/2,j}(u_{i+1,j}-u_{i,j}) - k_{i-1/2,j}(u_{i,j}-u_{i-1,j}) + k_{i,j+1/2}(u_{i,j+1}-u_{i,j}) - k_{i,j-1/2}(u_{i,j}-u_{i,j-1})\Big) + O(h).\end{aligned} \quad (70)$$

Under the assumption that $O(h)$ is negligibly small, we can then solve the above equation for $u_{i,j}$, obtaining

$$u_{i,j} \approx \frac{1}{K}\left(\frac{k_{i+1,j}}{k_{i+1,j}+k_{i,j}}u_{i+1,j} + \frac{k_{i-1,j}}{k_{i-1,j}+k_{i,j}}u_{i-1,j} + \frac{k_{i,j+1}}{k_{i,j+1}+k_{i,j}}u_{i,j+1} + \frac{k_{i,j-1}}{k_{i,j-1}+k_{i,j}}u_{i,j-1}\right), \quad (71)$$

where

$$K = \frac{k_{i+1,j}}{k_{i+1,j}+k_{i,j}} + \frac{k_{i-1,j}}{k_{i-1,j}+k_{i,j}} + \frac{k_{i,j+1}}{k_{i,j+1}+k_{i,j}} + \frac{k_{i,j-1}}{k_{i,j-1}+k_{i,j}}. \quad (72)$$

By constructing and solving Eq. (71) for all grid points, we can then compute the temperature at all grid points.

In the direct FDM analysis, we must solve a linear algebraic equation with a square matrix of order $n_{\text{FDM}}^G = \left(n_{\text{grid}}^G - 1\right)^2$ and a bandwidth of 4.

Here, we use the example of the equivalent conductivity for the middle grid point, $k$, in Eq. (71). In general, the thermal conductivity at the middle grid point cannot be defined correctly because the point is located at the interface between two small squares. In this case, the harmonic mean of the conductivities of the two squares can be regarded as the equivalent conductivity for the middle grid point.

We now consider the physical meaning of the harmonic mean using $k_{i+1/2,j}$ as an example. It is assumed that the temperature profile between the grid point $(i,j)$ and the middle grid point $(i + 1/2, j)$ and the temperature profile between the middle grid point $(i + 1/2, j)$ and the grid point $(i + 1, j)$ are both linear. The equivalent conductivity of the middle grid point $(i + 1/2, j)$, denoted by $k_{i+1/2,j}$, is then equal to the harmonic mean of the conductivities of the grid points $(i,j)$ and $(i + 1, j)$.

The above example can be illustrated using a numerical expression. From the equilibrium of the heat flow at the interface, where the middle grid point is $(i + 1/2, j)$, we obtain

$$k_{i,j}\frac{u_{i+1/2,j}-u_{i,j}}{0.5h} = k_{i+1,j}\frac{u_{i+1,j}-u_{i+1/2,j}}{0.5h} = k_{i+1/2,j}\frac{u_{i+1,j}-u_{i,j}}{h}. \quad (73)$$

When we solve this expression for $k_{i+1/2,j}$, we find

$$k_{i+1/2,j} = \frac{2k_{i+1,j}k_{i,j}}{k_{i+1,j}+k_{i,j}}. \quad (74)$$

We can then obtain the other three harmonic means in a similar manner.

$$k_{i-1/2,j} = \frac{2k_{i-1,j}k_{i,j}}{k_{i-1,j}+k_{i,j}},$$
$$k_{i,j+1/2} = \frac{2k_{i,j+1}k_{i,j}}{k_{i,j+1}+k_{i,j}}, \quad (75)$$
$$k_{i,j-1/2} = \frac{2k_{i,j-1}k_{i,j}}{k_{i,j-1}+k_{i,j}}.$$

## 5.3 SDM

We stated in Subsection 5.1 that $G$ consists of $\left(n_{\text{grid}}^G - 1\right) \times \left(n_{\text{grid}}^G - 1\right)$ small squares. Here, $n_{\text{grid}}^G$ is the number of grid points on a single side of the direct FDM model. For comparison with the direct FDM model and the SDM model, we set $n_{\text{grid}}^G$ as:

$$n_{\text{grid}}^G = 998. \quad (76)$$

Therefore, the total number of grid points in the direct FDM, denoted by $n_{\text{FDM}}^G$, is:

$$n_{\text{FDM}}^G = \left(n_{\text{grid}}^G - 1\right)^2 = 994009. \quad (77)$$

### 5.3.1 SDM global analysis

The global analysis is performed in accordance with the procedures that were explained in Subsection 3.2. We arrange the CPs at equal intervals of $h_{\text{SDM}}$ on $G$ (see Fig. 2(b)). As shown in Eq. (56), the total number of CPs in $G$, $n_{\text{SDM}}^G$, can be suppressed to be $(h/h_{\text{SDM}})^2$ times as many as $n_{\text{FDM}}^G$.

We tested the following three SDM models with different $h_{\text{SDM}}$ values, with these results:

SDM 1: $h_{\text{SDM}} = 4h$, $\quad n_{\text{SDM}}^G = 62500$
SDM 2: $h_{\text{SDM}} = 6h$, $\quad n_{\text{SDM}}^G = 27889$. $\quad (78)$
SDM 3: $h_{\text{SDM}} = 12h$, $\quad n_{\text{SDM}}^G = 7056$

The above indicates that all the $n_{\text{SDM}}^G$ values are considerably smaller than $n_{\text{FDM}}^G = 994009$.

The global influence coefficient matrix, $\mathbf{a}^G$, must then be established by assembling all the influence coefficient matrices, $\mathbf{a}^L$. We construct relational expressions for all CPs in $G$, in the form $\mathbf{u}^G = \mathbf{a}^G \mathbf{u}^G$, and solve these expressions based on the global boundary conditions. The order and the bandwidth of $\mathbf{a}^G$ are $n_{\text{SDM}}^G$ and 8, respectively.

### 5.3.2 SDM local analysis

To construct $\mathbf{a}^L$ and $\mathbf{N}^L$, a local analysis is conducted using the oversampling technique according to the procedures that were introduced in Subsection 3.1.2.

The following points correspond with those of the local analysis that was described in Subsection 4.3.2.
· The number of reference CPs is eight, as shown in Fig. 2(b).
· As in Fig. 4(b), each local domain $L^+$ is a square with a side length of $4h_{\text{SDM}}$ (i.e., it is four times as large as the interval between the CPs).
· The number of grid points in $L^+$ is $n_{\text{SDM}}^{L^+}$, as given by Eq. (59).

Conversely, the points below differ from the local analysis that was illustrated in Subsection 4.3.2.
· When analyzing $L^+$ using the FDM, we can use the same finite difference equation as that used for the direct FDM (Subsection 5.2 and Eq. (71)).
· The conductivity distribution is nonuniform throughout $G$ because the local domains are all different. There are

approximately $n_{\text{SDM}}^G$ different kinds of local domains in total.

To summarize the above, a linear algebraic equation that has a square matrix of order $n_{\text{SDM}}^{L+}$ must be solved approximately $n_{\text{SDM}}^G$ times when the local analysis is conducted.

## 5.4 DDM
### 5.4.1 Difference in example problem between DDM and SDM

We stated in Subsection 5.1 that $G$ is composed of $(n_{\text{grid}}^G - 1) \times (n_{\text{grid}}^G - 1)$ small squares. When the direct FDM model is compared with the DDM model, we set $n_{\text{grid}}^G$ to be
$$n_{\text{grid}}^G = 999. \tag{79}$$
Therefore, the number of small squares, $n_{\text{FDM}}^G$, is:
$$n_{\text{FDM}}^G = (n_{\text{grid}}^G - 1)^2 = 996004. \tag{80}$$
There is a small difference between the value of $n_{\text{FDM}}^G$ for the DDM analysis ($n_{\text{FDM}}^G = 996004$) and that of the SDM analysis ($n_{\text{FDM}}^G = 994009$, from Eq. (77)). Ideally, we should perform a complete analysis of the same $G$ using both the SDM and the DDM and then compare the results of the two methods. However, the $G$ of the DDM is slightly larger than the $G$ of the SDM. This is a result of the different ways in which $G$ is divided for the SDM and the DDM. Therefore, we cannot prepare an example problem in which $G$ has exactly the same structure.

However, the thermal conductivity distribution in the $G$ of the SDM, which consists of 997×997 squares, is exactly same as that in the (1st–997th rows)×(1st–997th columns) squares in the $G$ of the DDM. Therefore, these two global domains can be regarded as being almost the same.

### 5.4.2 Protocols

(1) First, divide $G$ into local domains that are squares with side lengths of $h_{\text{DDM}}$ in $L_{\text{DDM}}$, as depicted in Fig. 8.

(2) Then, construct fine-grained FDM models of $L_{\text{DDM}}$.
We use the same finite difference discretization here that was used for the direct FDM (see Subsection 5.2 and Eq. (71)).

(3) Set the outer and inner grid points.
As shown in Fig. 8(a), we define the grid points on the boundary of $L_{\text{DDM}}$ as the outer grid points (marked in yellow). In addition, the second outermost grid points are then defined as the inner grid points (marked in green). The temperatures at these grid points are then denoted by
$$\mathbf{u}_{\text{DDM}}^O = \left(u_1^O \cdots u_{n^O}^O\right)^T$$
$$\mathbf{u}_{\text{DDM}}^I = \left(u_1^I \cdots u_{n^I}^I\right)^T, \tag{81}$$

where $n^O$ is the number of outer grid points and $n^I$ is the number of inner grid points.
$$\begin{aligned} n^O &= 4\left(\frac{h_{\text{DDM}}}{h} - 1\right) \\ n^I &= 4\left(\frac{h_{\text{DDM}}}{h} - 2\right). \end{aligned} \tag{82}$$

(4) Conduct the local analysis of $L_{\text{DDM}}$.
The results of the FDM analysis of $L_{\text{DDM}}$ generate a matrix that determines the relationship between $\mathbf{u}_{\text{DDM}}^O$ and $\mathbf{u}_{\text{DDM}}^I$, denoted by $\mathbf{a}_{\text{DDM}}^L$. That is,
$$\mathbf{u}_{\text{DDM}}^I = \mathbf{a}_{\text{DDM}}^L \mathbf{u}_{\text{DDM}}^O. \tag{83}$$
$\mathbf{a}_{\text{DDM}}^L$ is an $n^I$-by-$n^O$ matrix. We conduct local analyses for all local domains in $G$ to obtain $\mathbf{a}_{\text{DDM}}^L$ for all these domains.

(5) Conduct the global analysis.
By constructing Eq. (83) for each of the local domains and then solving for each of them, we can compute the temperatures for all outer and inner grid points. This process is called the global analysis of the DDM.

### 5.4.3 Computational cost of the DDM
#### 5.4.3.1 Local analysis

The total number of grid points contained in each finite difference local domain is denoted by:
$$n_{\text{DDM}}^L = \left(\frac{h_{\text{DDM}}}{h} + 1\right)^2. \tag{84}$$
Additionally, the number of times that the local analysis was conducted (i.e., the number of different local domains in $G$) is denoted by:
$$m_{\text{DDM}}^L = \left(\frac{h}{h_{\text{DDM}} - h}\right)^2 (n_{\text{grid}}^G - 2)^2. \tag{85}$$
Thus, during the local analysis, we must solve a linear algebraic equation, which has a square matrix of order $n_{\text{DDM}}^L$, $m_{\text{DDM}}^L$ times.

#### 5.4.3.2 Global analysis

During the global analysis, the grid points in each local domain are deleted apart from the inner and outer grid points, as shown in the left figure in Fig. 8(a). The information contained on the deleted grid points is not used in the global analysis, and thus the computational cost is reduced by that amount. The total number of grid points in $G$ for the DDM is given by:
$$n_{\text{DDM}}^G = n_{\text{FDM}}^G - m_{\text{DDM}}^L \left(\frac{h_{\text{DDM}}}{h} - 3\right)^2. \tag{86}$$

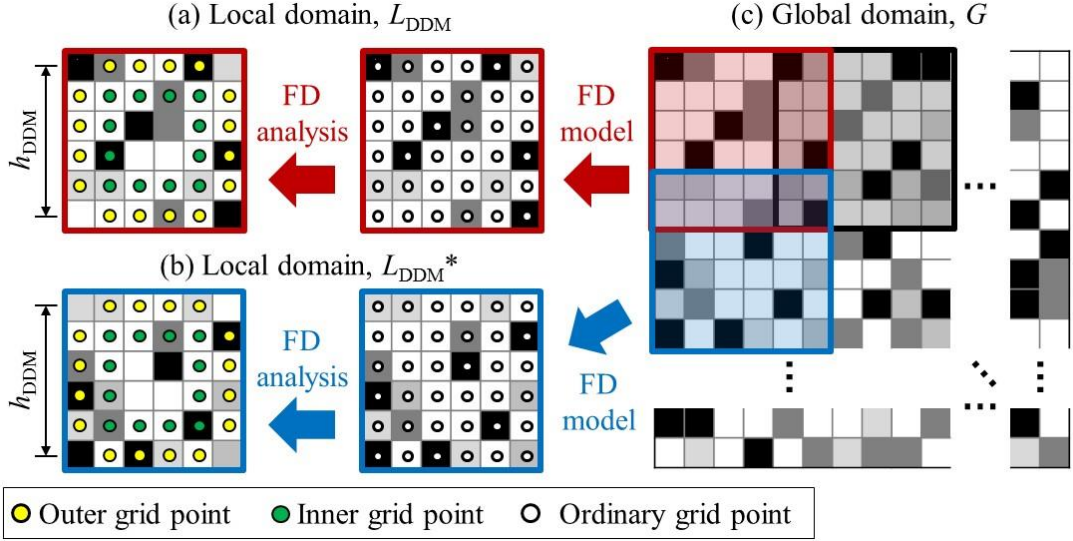

Fig. 8. (a) Local domain of the DDM, $L_{DDM}$; (b) another local domain adjacent to $L_{DDM}$; and (c) nonperiodic and heterogeneous global domain.

The equation above indicates that the number of deleted grid points increases as the side length of $L_{DDM}$ ($h_{DDM}$) increases, which then reduces the cost of the global analysis.

Conversely, the bandwidth of the square matrix in the global algebraic equation (and is referred to here as the global matrix) is approximately equivalent to the number of outer grid points of $L_{DDM}$, which is $n^O = 4(h_{DDM}/h - 1)$. The global matrix is thus strongly related to the scale of the global analysis. Therefore, an increase in $h_{DDM}$ increases both the bandwidth and the cost of the global analysis process.

We therefore tested the following five DDM models using different $h_{DDM}$ values.

DDM 1: $h_{DDM} = 4h$, $n^G_{DDM} = 748000$
DDM 2: $h_{DDM} = 6h$, $n^G_{DDM} = 555108$
DDM 3: $h_{DDM} = 12h$, $n^G_{DDM} = 307104$. (87)
DDM 4: $h_{DDM} = 83h$, $n^G_{DDM} = 51220$
DDM 5: $h_{DDM} = 166h$, $n^G_{DDM} = 27748$

The results above indicate that all values of $n^G_{DDM}$ are much smaller than $n^G_{FDM} = 996004$ in Eq. (80).

Comparison of Eqs. (78) and (87) shows that there is no difference between the numbers of divided regions of SDM $i$ and DDM $i$. Therefore,

$$h_{DDM} = h_{SDM}. \quad (88)$$

### 5.5 Results
#### 5.5.1 Calculation precision

The exact solution to this problem cannot be obtained. We therefore regard the RMSEs when compared with the direct FDM solution ($RMSE_{SDM}, RMSE_{DDM}$) as suitable accuracy indexes.

$$RMSE_{SDM} = \sqrt{\frac{1}{n^G_{SDM}} \sum_{i,j} \left(u_{SDM}(\mathbf{x}_{i,j}) - u_{FDM}(\mathbf{x}_{i,j})\right)^2}$$

$$RMSE_{DDM} = \sqrt{\frac{1}{n^G_{DDM}} \sum_{i,j} \left(u_{DDM}(\mathbf{x}_{i,j}) - u_{FDM}(\mathbf{x}_{i,j})\right)^2}. \quad (89)$$

The precision of the SDM is regarded as being higher as $RMSE_{SDM}$ becomes closer to 0.

$RMSE_{SDM} = 0$ only when $\mathbf{a}^L_{ERROR}$ in Eq. (27) is a zero matrix for all local domains. As stated in Subsection 2.2.2.3, the precision of the SDM does not exceed that of the direct FDM; therefore, $RMSE_{SDM} > 0$ in general. If $RMSE_{DDM} = 0$, then it can be said that the DDM solution corresponds exactly with the direct FDM solution.

$RMSE_{SDM}$ and $RMSE_{DDM}$ are shown in Fig. 9. For all the results gathered here,

$$RMSE_{SDM}, RMSE_{DDM} < 0.0002. \quad (90)$$

The minimum and maximum temperatures in the global field are 0 and 1°C, respectively. When compared with the temperature difference of 1°C, the errors above are sufficiently small.

SDM $i$ and DDM $i$ have local domains of the same size (i.e., $h_{SDM} = h_{DDM}$) and both of these domains become larger with increasing $i$. $RMSE_{SDM}$ and $RMSE_{DDM}$ remain almost constant, regardless of the value of $i$. However, there is a common factor in that SDM 2 generates the minimum error among all the SDM models and DDM 2 provides the minimum error among all the DDM models.

#### 5.5.2 Computation time

Figure 10 shows the analytical times (i.e., the CPU times) for the direct FDMs, SDMs, and DDMs. We show the CPU times that were consumed by the global analyses in gray and the times for the local analyses in black.

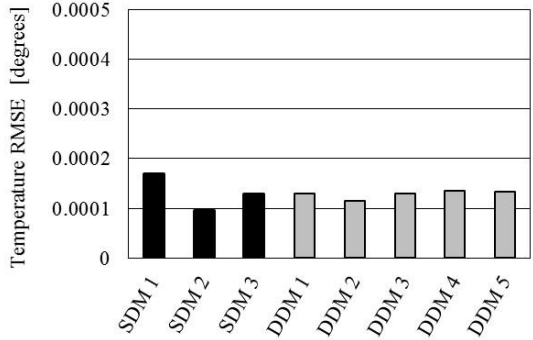

Fig. 9. Temperature RMSEs of the SDM models and the DDM models obtained from analysis of the nonperiodic stationary temperature field shown in Fig. 7.

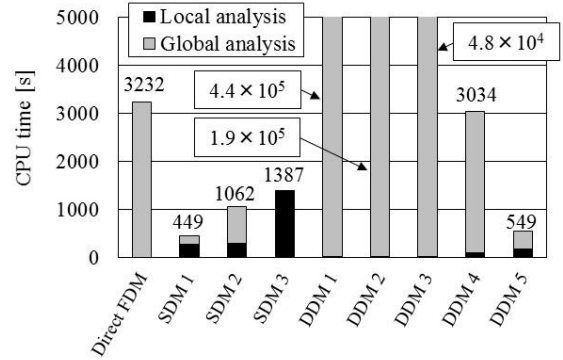

Fig. 10. CPU times consumed for SDM analyses and DDM analyses of the nonperiodic stationary temperature field shown in Fig. 7.

When compared with the total time for the direct FDM, DDM 4 is equivalent; however, SDMs 1–3 and DDM 5 are computationally lower in cost. In particular, SDM 1 and DDM 5 reduce the total time to approximately one-sixth of that for the direct FDM. Conversely, DDMs 1–3 require enormous amounts of time that significantly exceed that required for the direct FDM.

We must therefore discuss the reasons why the CPU times above were obtained.

When the sizes of the local domains ($h_{SDM}$ and $h_{DDM}$) increase, the changes in the costs for both the SDM and the DDM are as follows:
· the number of degrees of freedom in $G$ decreases, which then reduces the cost of the global analysis;
· the number of local domain types present decreases, which thus reduces the number of times that the local analysis must be conducted;
· the number of grid points in the local domains increases with any increase in the local domain size, which then increases the cost per performance for the local analysis process.

The CPU times per single local analysis process for each of the methods are compared in Fig. 11. When compared with the local domain of the DDM, the local domain of the SDM when using the oversampling method is four times longer and has 16 times (or 64 times in the 3D field case) the number of grid points. Thus, the cost of the local analysis for SDM $i$ increases rapidly with increasing $i$.

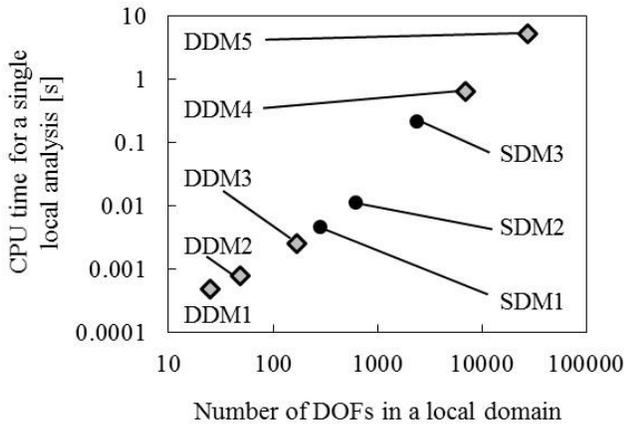

Fig. 11. Comparison of CPU times taken for each local analysis between the SDM and the DDM.

Figure 10 demonstrates that the time taken for the local analysis is much longer than that for the global analysis in SDM 3. SDMs 4 and 5 cannot be tested in this study because of the enormous amount of time required for the local analyses.

In each of DDMs 1–5, the global analysis time is longer than that required for the local analysis. In particular, DDMs 1–3 are vastly inferior to the direct FDM in terms of total CPU time. This is because DDMs 1–3 require only a short time for the local analysis but a much longer time for the global analysis. When compared with the direct FDM, the DDM's global domain has fewer degrees of freedom, but the DDM's global analysis requires greater numbers of temperature values to be provided at the grid points for reference during grid point temperature estimation. In DDMs 1–3, the cost increases related to increasing the bandwidth of the global equation would be much greater than the cost savings related to reduction of the degrees of freedom.

To minimize the total time for both the SDM and the DDM, we must determine the number of divided regions for the global domain that would produce approximately equivalent CPU times for both the local and global analyses.

## 6 Conclusions

In this work, we have provided a detailed mathematical analysis of the error factors of the multiscale SDM, which have not been revealed previously.

We investigated the convergence properties of the SDM using a linear stationary temperature problem and compared the results with those obtained when using a standard FDM and a conventional DDM [7,8].

Additionally, we applied the SDM to the analysis of a linear temperature field with a nonuniform and nonperiodic thermal conductivity field for the first time. We subsequently solved the same type of problem using the DDM. We then compared the performances of the SDM and the DDM. The results were as follows.

- The accuracy of the SDM is very high, and is approximately equivalent to that of the DDM (the RMSE for the temperature is less than 0.02% of the maximum temperature in the simulated field).
- From the calculated results for three SDM models, the total CPU time is 13–43% of that of the direct FDM model.
- The most efficient DDM model requires 17% of the CPU time of the direct FDM model.
- To minimize the total computation time required for both the SDM and DDM, we must determine the number of divided regions for the global domain that provides approximately equivalent CPU times for the local and global analyses.

## Acknowledgments

Funding: This work was supported by the Japan Society for the Promotion of Science (JSPS) KAKENHI [grant number 17K14144].

## Nomenclature

| Symbol | Description |
|---|---|
| $\mathbf{a}^G$ | Influence coefficient matrix for global domain of seamless-domain method (SDM), $G$ |
| $\mathbf{a}^L$ | Influence coefficient matrix for SDM's local domain |
| $d \in \{1, \ldots, 3\}$ | Number of dimensions of a domain |
| $f(\mathbf{x})$ | Source value at point $\mathbf{x}$ |
| $G \subset \mathbf{R}^d$ | Global domain |
| $h_{\text{SDM}}$ | Interval between coarse-grained points (CPs) in SDM model |
| $h_{\text{DDM}}$ | Length of local domain in domain decomposition method (DDM) model |
| $h$ | Interval between grid points in finite-difference method (FDM) model |
| $L \subset G$ | SDM's local domain without oversampled domain |
| $L^+ \supset L$ | SDM's local domain including oversampled domain |
| $L^* \subset \mathbf{R}^d$ | SDM's local domain next to $L$ |
| $L_{\text{DDM}} \subset G$ | DDM's local domain |
| $m_{\text{DDM}}^L$ | Number of local domains in DDM's global domain |
| $n_{\text{CP}}^G$ | Number of CPs on one side of SDM's global domain |
| $n_{\text{DDM}}^G$ | Number of grid points in DDM's global domain |
| $n_{\text{FDM}}^G$ | Number of grid points in direct FDM's global domain |
| $n_{\text{grid}}^G$ | Number of grid points on one side of direct FDM's global domain |
| $n_{\text{SDM}}^G$ | Number of CPs in SDM's global domain |
| $n_{\text{DDM}}^L$ | Number of grid points in DDM's local domain, $L_{\text{DDM}}$ |
| $n_{\text{SDM}}^{L+}$ | Number of grid points in SDM's local domain, $L^+$ |
| $\mathbf{N}^L$ | Interpolation function matrix for dependent variable in $L$ |
| $\mathbf{N}^{\partial L}$ | Interpolation function matrix for dependent variable on $L$'s boundary, $\partial L$ |
| **R** | Set of all real numbers |
| $RMSE_{\text{DDM}}$ | Root mean squared error of DDM solution |
| $RMSE_{\text{FDM}}$ | Root mean squared error of FDM solution |
| $RMSE_{\text{SDM}}$ | Root mean squared error of SDM solution |
| $u(\mathbf{x})$ | Dependent-variable value at point $\mathbf{x}$ |
| $u_{\text{CP}i} = u(\mathbf{x}_{\text{CP}i})$ | Dependent variable for the $i$th CP |
| $\mathbf{u}^G$ | Dependent variable for all CPs in global domain, $G$ |
| $u_{i,j} = u(\mathbf{x}_{i,j})$ | Dependent variable for grid point $(i,j)$ |
| $\mathbf{u}^L$ | Dependent variable for all CPs on $L$'s boundary, $\partial L$ |
| $\mathbf{u}^{L^+}$ | Dependent variable for all boundary points on $L^+$'s boundary, $\partial L^+$ |
| $\mathbf{u}^I_{\text{DDM}}$ | Dependent variable for inner grid points in DDM's local domain |
| $\mathbf{u}^O_{\text{DDM}}$ | Dependent variable for outer grid points in DDM's local domain |
| $\mathbf{x} \in \mathbf{R}^d$ | Position vector |
| $\mathbf{x}_{\text{CP}i}$ | Position of the $i$th CP |
| $\mathbf{x}_{i,j}$ | Position of grid point $(i,j)$ |
| $\partial G$ | Boundary of $G$ |
| $\partial L$ | Boundary of $L$ |
| $\partial L^+$ | Boundary of $L^+$ |
| $\partial L^*$ | Boundary of $L^*$ |
| $\partial L_{\text{DDM}}$ | Boundary of $L_{\text{DDM}}$ |
| $\#^{-1}$ | Inverse of # |
| $\#^T$ | Transpose of # |